\documentclass[preprint,12pt]{elsarticle}
\usepackage{amssymb}
\usepackage{amsmath}
\newtheorem{theo}{Theorem}[section]

\newtheorem{cor}[theo]{Corollary}
\newtheorem{prop}[theo]{Proposition}
\newtheorem{pf}[theo]{Proof} 

\newtheorem{remark}[theo]{Remark}
\newcommand{\be}{\begin{eqnarray*}}
\newcommand{\ee}{\end{eqnarray*}}
\newcommand{\ben}{\begin{eqnarray}}
\newcommand{\een}{\end{eqnarray}}
\def\A {\mathbf{A}}

\def\I {\mathbf{I}}

\def\P {\mathbf{P}}

\def\R {\mathbf{R}}

\def\b {\boldsymbol{b}}

\def\l {\boldsymbol{l}}

\def\s {\boldsymbol{s}}

\def\u {\boldsymbol{u}}

\def\w {\boldsymbol{w}}
\def\x {\boldsymbol{x}}
\def\y {\boldsymbol{y}}

\def\Bc {\mathcal{B}}

\def\Nc {\mathcal{N}}

\def\Rb {\mathbb{R}}
\def\Pb {\mathbb{P}}

\begin{document}

\begin{frontmatter}

\author[els]{D.Ounaissi}
\ead{daoud.ounaissi@u-pec.fr}
\fntext[fn1]{Corresponding Author}
\author[els]{N.Rahmania}
\ead{nadji.rahmania@univ-lille1.fr}
\address{Laboratoire LAMA, UFR de Math\'ematiques, Universit\'e de Paris 12, France.\\ Laboratoire Paul Painlev\'e, UFR de Math\'ematiques, Universit\'e de Lille, France.}

\title{Bayesian Lasso : Concentration and MCMC Diagnosis}

\begin{abstract}
Using posterior distribution of Bayesian LASSO  
we construct a semi-norm on the parameter space.  
We show that the partition function depends on the ratio
of the $l^1$ and $l^2$ norms and present three regimes.  
We derive the concentration of Bayesian LASSO, and 
present MCMC convergence diagnosis. 

\end{abstract}

\end{frontmatter}

\section{Introduction} 
Let $p \geq n$ be two positive integers, $\y \in \Rb^n$ and 
$\A$ be an $n\times p$ matrix with real numbers entries. Bayesian LASSO 
\ben 
c(\x)=\frac{1}{Z}\exp\Big(-\frac{\|\A\x-\y\|_2^2}{2}-\|\x\|_1\Big)\label{c}
\een   
is a typically posterior distribution used in the linear regression 
\be 
\y=\A\x+\w. 
\ee 
Here 
\ben 
Z=\int_{\Rb^p}\exp\Big(-\frac{\|\A\x-\y\|_2^2}{2}-\|\x\|_1\Big)d\x\label{Zx}
\een

is the partition function, $\|\cdot\|_2$ and $\|\cdot\|_1$ are respectively 
the Euclidean and the $l_1$ norms. 
The vector $\y\in\Rb^n$ are the observations, $\x\in \Rb^p$ is the unknown signal to recover, $\w\in\Rb^n$ is the standard Gaussian noise, and $\A$ is a known matrix which maps the signal domain $\Rb^p$ into the observation domain $\Rb^n$. If we suppose that $\x$ is drawn from 
Laplace distribution i.e. the distribution proportional to 
\ben 
\exp(-\|\x\|_1),
\label{laplace}  
\een
then the posterior of $\x$ known $y$ is drawn from the distribution $c$ (\ref{c}). 
The mode  
\ben 
\arg\min\Big\{\frac{\|\A\x-\y\|_2^2}{2}+\|\x\|_1:\quad \x\in\Rb^p\Big\}
\een
of $c$ was first introduced in \cite{Tibshirani1996} and called LASSO. It is also called Basis Pursuit De-Noising method \cite{Chen}. In our work we select the term LASSO and keep it for the rest of the article.  

In general LASSO is not a singleton, i.e. the mode of the distribution $c$ 
is not unique. In this case LASSO is a set and we will denote by lasso any element of this set. 
A large number of theoretical results has been provided 
for LASSO. See \cite{Daubechies2004}, \cite{DDN}, \cite{Fort}, \cite{Mendoza}, \cite{Pereyra} and the references herein. 
The most popular algorithms to find LASSO are LARS algorithm \cite{Efron}, ISTA and FISTA algorithms see e.g. \cite{Beck} and the review article \cite{Parikh}.

The aim of this work is to study geometry of bayesian LASSO 
and to derive MCMC convergence diagnosis. 
\section{Polar integration}   
Using polar coordinates $\s=\frac{\x}{\|\x\|}\in\mathcal{S}, r=\|\x\|$, the partition function (\ref{Zx})  
\ben 
Z_p=\int_{S}Z_p(\s)d\s,
\label{Zp} 
\een 
where $\|\cdot\|$ denotes one of  $l^2$ or $l^1$  norms
in $\Rb^p$, $d\s$ denotes the surface measure on the unit sphere $S$ of the norm $\|\cdot\|$, and 
\ben 
Z_p(\s)=\int_0^{+\infty}\exp\{-f(r\s)\}r^{p-1}dr, 
\label{Zs}
\een 
where $f(\x):=\frac{\|\A\x-\y\|_2^2}{2}+\|\x\|_1, \x\in\Rb^p.$

We express the partition function (\ref{Zs}) using the parabolic cylinder function. We also give an inequality of concentration and a geometric interpretation of the partition function $ Z_p $.
 
\section{Parabolic cylinder function and partition function} 
We extend the function $\s\in \mathcal{S}\to Z_p(\s)$ à 
\ben\label{Zpx} 
\x\in\Rb^p\to Z_p(\x)=\int_0^{+\infty}\exp\{-f(r\x)\}r^{p-1}dr.
\een 
This extension is homogeneous of order $ -p $.

If  $\A\x=0$, then $f(r\x)=\frac{\|\y\|_2^2}{2}+r\|\x\|_1$, and more if $\x\neq 0$, then
\be 
Z_p(\x)=(p-1)!\|\x\|_1^{-p}\exp(-\frac{\|\y\|_2^2}{2}).
\ee 
If $ \A \x \neq 0 $, then we will express $ Z_p (\x) $ using the parabolic cylinder function. We recall that for $ \b \in \R , z \geq 0$ the parabolic cylinder function is given by
\ben\label{Dasymptoticylindre}   
D_b(z)=z^{b}\exp(-\frac{z^2}{4})[1+O(z^{-2})],
\een 
when $z\to +\infty$ \cite{Temme}. We also recall the integral representation of Erdélyi \cite{Erdélyi} for the parabolic cylinder function
\be 
\exp(\frac{z^2}{4})\Gamma(\nu) D_{-\nu}(z)=\int_0^{+\infty}\exp(-\frac{1}{2}r^2-zr)r^{\nu-1}dr,
\quad \nu >0,
\ee 
where $\Gamma(\nu)=\int_0^{+\infty}\exp(-t)t^{\nu-1}dt$  
is the $\Gamma$ function. 
\begin{prop} 
The variable
\ben\label{omegalasso}  
\omega_{lasso}:=\frac{\|\x\|_1-\langle\A\x,\y\rangle}{\|\A\x\|_2}
\een 
will play an important role. It depends only on $\s=\frac{\x}{\|\x\|_1}\in\mathcal{S}_{p-1,1}$ and the function $\s\in\mathcal{S}_{p-1,1}\to \omega_{lasso}(\s) $ is bounded below by $\lambda_{1,2}:=\min\{\frac{1}{\|\A\s\|_2}-\frac{\langle\A\s,\y\rangle}{\|\A\s\|_2}:\quad \s\in 
\mathcal{S}_{p-1,1}\}. $
\end{prop}

Now we can announce the following result.  
\begin{prop}
We have for $\A\x \neq 0$
\be 
Z_p(\x)=(p-1)!\exp(-\frac{\|\y\|_2^2}{2})\|\A\x\|_2^{-p}
\exp(\frac{\omega_{lasso}^2}{4})D_{-p}(\omega_{lasso}).
\ee 
If $\A\x\to 0$, then $\omega_{lasso}\to +\infty$ and  
\be  
Z_p(\x)=(p-1)!\exp(-\frac{\|\y\|_2^2}{2})
\|\x\|_1^{-p}
[1+O(\omega_{lasso}^{-2})]. 
\ee 
\end{prop} 

\begin{cor}
If $y=0$ then $\omega_{lasso}=\frac{1}{\|\A\s\|_2}$ is bounded below by $\frac{1}{\lambda_{1,2}}$, where
$\lambda_{1,2}=\max(\|\A\s\|_2:\s\in\mathcal{S}_{p-1,1})$ is the norm of the operator $\A:(\Rb^p,\|\cdot\|_1)\to 
(\Rb^n,\|\cdot\|_2)$. The partition fucntion
\be 
Z_p(\s)=(p-1)!\omega_{lasso}^{p}
\exp(\frac{\omega_{lasso}^2}{4})D_{-p}(\omega_{lasso})
\ee 
is $\|\A\s\|_2^2$ convex and decreasing.  
\end{cor} 
\begin{pf}
It suffices to remark that $Z_p(\s)=\int_0^{+\infty}\exp\{-\frac{\|\A\s\|^2}{2}r^2-r\}r^{p-1}dr. $
\end{pf}
\section{Geometric interpretation of the partition function} 
First we represent $f(r\x)$ for $\A\x\neq 0$ in the form
\ben\label{frbeta} 
f(r\x)=\frac{\|\y\|_2^2}{2}-\frac{\omega_{lasso}^2}{2}+\frac{(r\|A\x\|_2+\omega_{lasso})^2}{2}. 
\een 
The function $\exp\{-f(\x)\}, \forall \x\in\Rb^p$ is log-concav and integrable in $\Rb^p$.  Observe that $Z_p^{-\frac{1}{p}}$ is a norm on the null space $Ker(\A)$ of $\A$. A general result \cite{Ball} tells us that
\be 
\x\in \Rb^p\to Z_p^{-\frac{1}{p}}(\x):=\|\x\|_c
\ee
is a quasi-norm on $\R^p$.  The unit ball of $ \| \cdot \|_c $ is defined by
\be 
\mathcal{B}(\A,\y)&:=&\{\x\in\Rb^p: \|\x\|_c\leq 1\}\\
&=&\{\x\in\Rb^p: Z_p(\x)\geq 1\}\\
&=&\{\x=r\s \in\Rb^p: Z_p(\s)\geq r^p\}\\
&=&\{\x=r\s \in\Rb^p:\,(p-1)! \exp(-\frac{\|\y\|_2^2}{2})\|\A\s\|_2^{-p}\exp(\frac{\omega_{lasso}^2}{4})D_{-p}(\omega_{lasso})\geq r^p\}. 
\ee 
Its contour is equal to
\be 
C(\A,\y)&:=&\{\x\in\Rb^p: \|\x\|_c=1\}\\
&=&\{\x\in\Rb^p: Z_p(\x)=1\}\\
&=&\{\x=r\s\in\Rb^p: Z_p(\s)= r^p\}\\
&=&\{\x=r\s\in\Rb^p:\,(p-1)! \exp(-\frac{\|\y\|_2^2}{2})\|\A\s\|_2^{-p}\exp(\frac{\omega_{lasso}^2}{4})D_{-p}(\omega_{lasso})= r^p\}. 
\ee 
We summarize our results in the following proposition.
\begin{prop}
1) For each $\s \in \mathcal{S}_{p-1,1}$, the longest segment $ [0, r]\s $ contained in $ \mathcal{B}(\A, \y) $ holds for $ r = r_{max} (\s) $ is solution of the equation
\be 
r^p=(p-1)!\exp(-\frac{\|\y\|_2^2}{2})\|\A\s\|_2^{-p}\exp(\frac{\omega_{lasso}^2}{4})D_{-p}(\omega_{lasso}).
\ee 
2) The ball 
\be
\mathcal{B}(\A,\y)=
\bigcup_{\s \in \mathcal{S}_{p-1,1}}[0,r_{max}(\s)]\s,
\ee
and its contour is equal to 
\be 
C(\A,\y)=\{r_{max}(\s)\s:\quad  \s\in \mathcal{S}_{p-1,1}\}.
\ee 
3) The volume  $\mathcal{B}(\A,\y)$ is
\be 
\int_{\mathcal{S}_{p-1,1}} \frac{r_{max}^p(\s)}{p}d\s=\frac{Z_p}{p}. 
\ee
\end{prop}

\section{Necessary and sufficient condition to have lasso equal zero}
Now we can give the necessary and sufficient condition to have\\ $lasso=\{0\}$
\begin{prop}
The following assertions are equivalent.

1) $0=\mbox{lasso}$.

2) $\omega_{lasso}(\s)\geq 0$ pour tout $\s\in \mathcal{S}_{p-1,1}$. 

3) $\|\A^{\top}\y\|_{\infty} \leq 1.$
\end{prop}

\section{Concentration around the lasso}
\subsection{The case lasso null} 
The polar coordinate formula tells us that, we can draw a vector $ \x $ from $ c(\x)d\x$ by drawing its angle $ \s $ uniformly, and then simulate its distance $r$ to the origin from
\ben\label{crs} 
c(r,\s)dr=\frac{1}{Z_p(\s)}\exp\{-f(r\s)\}r^{p-1}dr 
\een 
Now let's estimate for $ r> 0 $ the probability
\be 
\P(\|\x\|> r)=\int_{\mathcal{S}}\int_r^{+\infty} 
c(r,\s)dr\frac{d\s}{|\mathcal{S}|},
\ee 
where $|\mathcal{S}|$ denotes the Lebesgue measure of $\mathcal{S}$. 
We introduce for each pair $ a \geq 0 $, $ b \in \Rb^p $ the function
\ben\label{gabp} 
g_{a,b,p}(r):=g_{a,b}(r)-(p-1)\ln(r),\quad r>0. 
\een
In the following
\be 
a=\|\A\s\|_2,\quad b=\omega_{lasso}.
\ee 

The function $r\geq 0\to g_{a,b}(r)$ is increasing (because
$b:=\omega_{lasso} \geq 0$). 
The function
$r\to g_{a,b,p}(r)$ 
f est convexe et atteint son minimum au point $r(\s)$ solution  de l'équation 
\be 
\|\A\s\|_2(r\|\A\s\|_2+\omega_{lasso})=\frac{p-1}{r}.
\ee 
The positive root is given by
\ben 
r(\s)=\frac{-\omega_{lasso}+\sqrt{\omega_{lasso}^2+4(p-1)}}{2\|\A\s\|_2}
\label{romega} 
\een
On one hand
\be 
\int_0^{+\infty}\exp\{-g_{a,b,p}(r)\}dr&\geq &
\exp\{-g_{a,b}(r(\s))\} 
\int_0^{r(\s)}r^{p-1}dr=\exp\{-g_{a,b,p}
(r(s))\}\frac{r(\s)}{p}. 
\ee 
On the other hand by using the convexity of
$r\to g_{a,b}(r)$, we have for all $r>0$, 
\be 
g_{a,b}(r)\geq g_{a,b}(r(\s))+\frac{(p-1)(r-r(\s))}{r(\s)},
\ee  
because $\partial_rg_{a,b}(r(\s))=\frac{p-1}{r(\s)}$. 
We deduce for $q>0$,     
\be 
\int_{qr(\s)}^{+\infty}\exp\{-g_{a,b,p}(r)\}dr 
\leq \exp\{-g_{a,b}(r)+p-1\}\int_{qr(\s)}^{+\infty} 
\exp\{-\frac{p-1}{r(\s)}r\} r^{p-1}dr\\
\leq \exp\{-g_{a,b}(r(\s))+p-1\}\int_{q(p-1)}^{+\infty} 
\exp(-r)r^{p-1}dr \frac{r(\s)}{(p-1)^p}\\
\leq \exp\{-g_{a,b}(r(\s))+p-1\}\frac{r(\s)}{(p-1)^p}\Gamma(p,q(p-1)), 
\ee 
where $\Gamma(\nu,x)=\int_x^{+\infty}\exp(-t)t^{\nu-1}dt$ is the upper incomplete gamma function. Finally we get the following result. 

\begin{prop}\label{0concentration} We have for all $q >0$,  
\ben\label{Pqp} 
\P(\|\x\| \geq qr(\s))\leq \frac{p\exp(p-1)}{(p-1)^p}\Gamma(p,q(p-1)):=P(q,p). 
\een 
\end{prop} 
Using the following estimate \cite{Natalini}
\be 
x^{a-1}\exp(-x) < \Gamma(a,x) < B x^{a-1} \exp(-x),\quad \forall\, a >1,\,\,B >1,\,\,x > \frac{B}{B-1}(a-1),
\ee 
we get for $q >1$, 
\be 
\Gamma(p,q(p-1))\leq 2 q^{p-1}(p-1)^{p-1}\exp(-q(p-1)).
\ee 
Therefore the quantity 
\be 
P(q,p)\leq \frac{2pq^{p-1}}{(p-1)}\exp\{-(q-1)(p-1)\}.
\ee

{\bf balance sheet .} If $\x$ is drawn from the density $c$, alors 
$\frac{\x}{r(\theta)}\in \mathcal{B}_2(0,q)$ with a probability at least equal to
$1-P(q,p)$. 

In the figure(1) we plot for  $p=2$, $n=1$, $\A= \begin{pmatrix} 1&1\end{pmatrix}$ and $\y=0$ the density $c(r,\s)dr$ for a fixed value of $\s$.\\

\begin{figure}[!ht]\label{crsplot}
\centering\includegraphics[width=14cm]{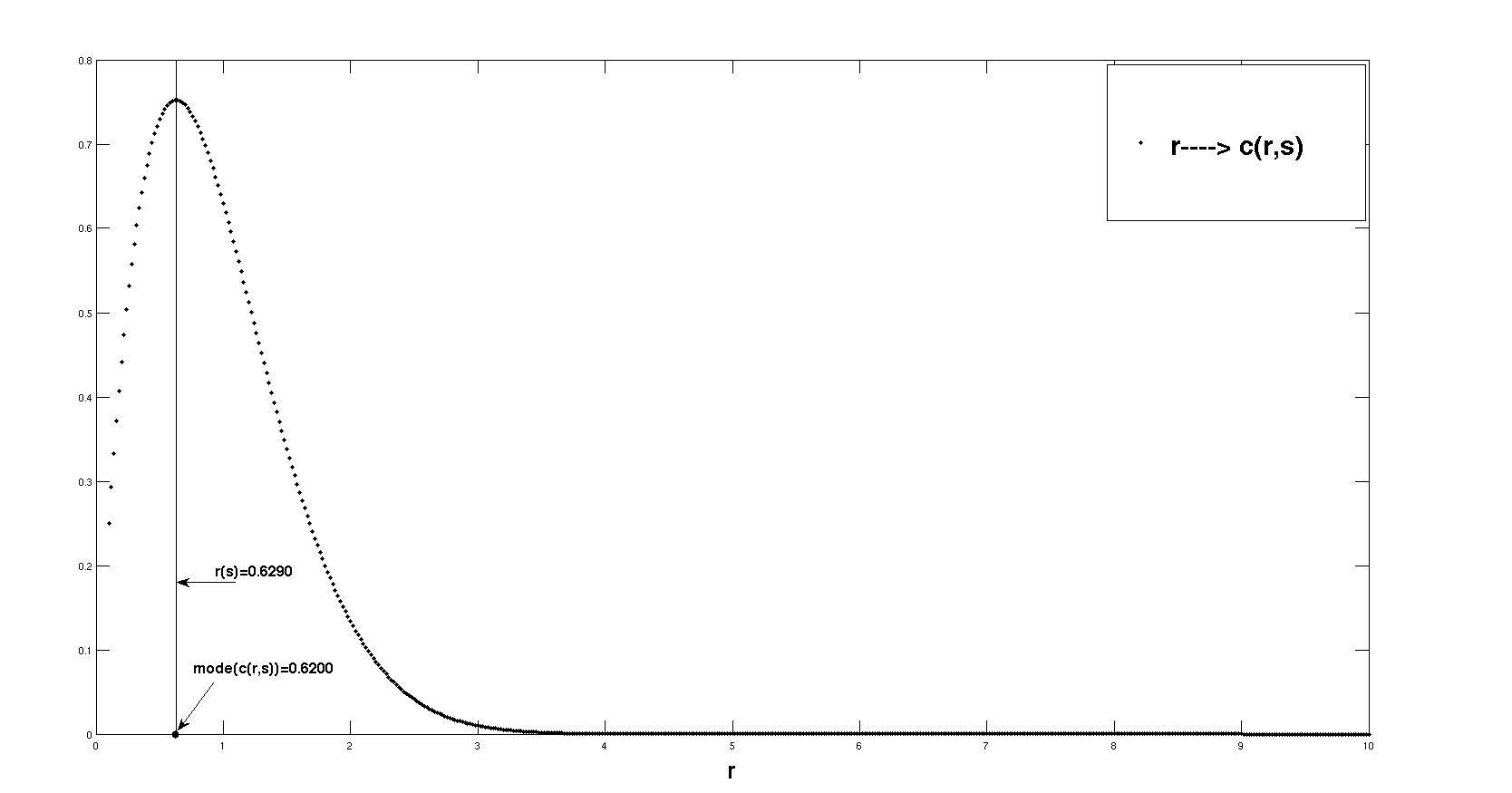} 
\caption{pour $r \in [0.1; 10] \to c(r, \s)$. }
 \end{figure} 
 
We notice that the mode $ c (r, \ s) = 0.6200$ is very close to the value of $ r (\ s)=0.6290$ (\ref {romega}) for the same fixed $ \s$.
\subsection{The general case}
We take the vector $\l\in lasso$. 
We will study the concentration of  $c$ around $\l$. 
The variable of interest is $\u = \x- \l $. The law of $\u $ has for density
\be 
c(\u+\l)d\u=\frac{1}{Z_p}\exp\{-f(\u+\l,\A,\y)\}d\u. 
\ee 
The change of variables formula gives for each norm $\|\cdot\|$ 
\be 
c(\u+\l)d\u=\frac{1}{Z_p}\exp\{-f(r\theta+\l)\}r^{p-1}drd\theta,\quad r >0,\, \theta\in \mathcal{S}. 
\ee 
By definition for any vector $\x$, the convex function $r \geq 0 \to f (r\s+ \l)$ reaches its minimum at the point $r=0$. Therefore
$r\geq 0\to f(r\s+\l)$ is increasing.  

The function  
\ben\label{frthetap} 
f(r\s+\l,p):=f(r\s+\l)-(p-1)\ln(r),\quad r >0, 
\een 
is strictly convex. Its critical point $r_{\l}(\s)$ is solution of the equation
\be 
\partial_rf(r\s+\l)=\frac{p-1}{r}.
\ee 
By a similar proof to that of propostion (\ref{0concentration}) we have the following result; 

\begin{prop} If $\x$ is drawn from the density $c$,
and $\s=\frac{\x-\l}{\|\x-\l\|}$, then for all $q >0$,  
\ben\label{Pqp} 
\P(\|\x-\l\|\geq qr_{\l}(\s))\leq \frac{p\exp(p-1)}{(p-1)^p}\Gamma(p,q(p-1)):=P(q,p). 
\een 
\end{prop} 
\section{Applications} 
\subsection{The contour in the case $p=2$, $n=1$}
Let $\A:=(a_1,a_2)$ a matrix of order $1\times 2$. Its null-space $Ker(\A)=\{(x_1,x_2):\quad a_1x_1+a_2x_2=0\}$.  We have that $\mathcal{B}(a_1,a_2,y)$ contains
\be 
Ker(\A)\bigcap \mathcal{B}_{2,1}. 
\ee 
This intersection is a symmetric segment noted $[(x_1(a_1,a_2),x_2(a_1,a_2)), -(x_1(a_1,a_2),x_2(a_1,a_2))]$.  

To determine the other points of the set $\mathcal{B}(a_1,a_2,y)$, 
we will directly calculate  $Z_2(\s)$. A simple calculation gives 
\be
Z_2(\s)=\exp(-\frac{y^2}{2}+\frac{\omega_{lasso}^2}{2})\int_0^{+\infty}
\exp\{-\frac{(|\A\s| r+\omega_{lasso})^2}{2}\} rdr,
\ee 
et 
\be 
|\A\s|\int_0^{+\infty}\exp\{-\frac{(|\A\s| r+\omega_{lasso})^2}{2}\} rdr+\omega_{lasso}\int_0^{+\infty}\exp\{-\frac{(|\A\s| r+\omega_{lasso})^2}{2}\}dr=1.
\ee 
Finally we have the following proposition. 
\begin{prop}\label{p=2} 1) If $\A\s\neq 0$, then 
\be 
Z_2(\s)=\exp(-\frac{y^2}{2}+\frac{\omega_{lasso}^2}{2})|\A\s|^{-1}
\{1-\frac{\omega_{lasso}}{|\A\omega|}\sqrt{2\pi}(1-F(\omega_{lasso}))\}, 
\ee 
where $F$ is the distibution function of the normal law. 

2) If $\A\s\neq 0$ and $y=0$, then 
\be 
Z_2(\s)&=&\omega_{lasso}\exp(\frac{\omega_{lasso}^2}{2})\{1-\omega_{lasso}^2\sqrt{2\pi}(1-F(\omega_{lasso}))\}.
\ee 

3) Ifi $\s\in \mathcal{S}_{1,1}$, $\A\s\neq 0$ and $y=0$, then
the function $z_2$
\be 
z_2(b^2)=\frac{1}{b}\exp(\frac{1}{2b^2})\{1-\frac{1}{b^2}\sqrt{2\pi}
(1-F(\frac{1}{b}))\}
\ee 
defined on $(0, \lambda_{1,2}^2]$ is convex and decreasing,
where $\lambda_{1,2}=\max_{\s\in \mathcal{S}_{1,1}}|\A\s|$. 

4) We have for $\s \in \mathcal{S}_{1,1}$  
\be 
Z_2(\s)=z_2(\frac{1}{\omega_{lasso}^2}),\quad \forall\,\omega\in \mathcal{S}_{1,1}.
\ee 
\end{prop} 

the ball 
\be 
\mathcal{B}(\A,0)=\{r\s: \quad Z_2(\s)\geq r^2\}. 
\ee 
is contained in the unit disk  $\|\x\|_1\leq 1$ for the norm $l^1$.  
The contour is defined by the equation 
\be 
Z_2(\s)=r^2.
\ee 
The norm of the linear operator $\A:(\Rb^2, \|\cdot\|_1)\to (\Rb, \|\cdot\|_2)$ 
is defined by 
\be 
\lambda_{1,2}=\sup_{\s:\|\s\|_1=1}\|\A\s\|_2. 
\ee 
the function $\s\to Z_2(\s)=z_2(\lambda_{1,2}^2)$ is constant on est constante sur  
\be 
\Omega_{1,2}=\{\s:\quad \|\s\|_1=1,\quad \|\A\s\|=\lambda_{1,2}\}. 
\ee 
If $\A=(1,1)$ then 
\be 
\Omega_{1,2}&=&\{\s:\quad \|\s\|_1=1,\quad \|\A\s\|=\lambda_{1,2}\}\\
&=&[(1,0), (0,1)]\bigcup [(-1,0), (0,-1)].
\ee 
If $\A=(a_1,a_2)$ with $|a_1| < |a_2|$, then 
\be 
\Omega_{1,2}&=&\{\s:\quad \|\s\|_1=1,\quad \|\A\s\|_2=\lambda_{1,2}\}\\
&=&\{(0,sgn(a_2)), (0,-sgn(a_2))\}. 
\ee 
In both case
\be 
\{z_2(\lambda_{1,2}^2)\}^{\frac{1}{2}}\Omega_{1,2} 
\ee   
is part of the contour. The other points of the contour are deduced from the equation
\be 
z_2(b^2)=a^2, \quad b\in (0, \lambda_{1,2}).  
\ee 
Each pair $(a,b)$ generate four points of $\mathcal{B}((a_1,a_2),0)$ of the form
$a\s$ where 
\be 
|\s_1|+|\s_2|=1,\quad 
|a_1\s_1+a_2\s_2|=b. 
\ee 

We plot in the figure 2 the contour of $\mathcal{B}(a_1,a_2,0)$ for different choices of the matrix $(a_1, a_2)$. We notice that the surface of $\mathcal{B}((a_1,a_2),0)$ is decreasing function relatively the norm $\lambda_{1,2}$ of the matrix $\A$.
\begin{figure}[!ht]\centering
\includegraphics[width=14cm]{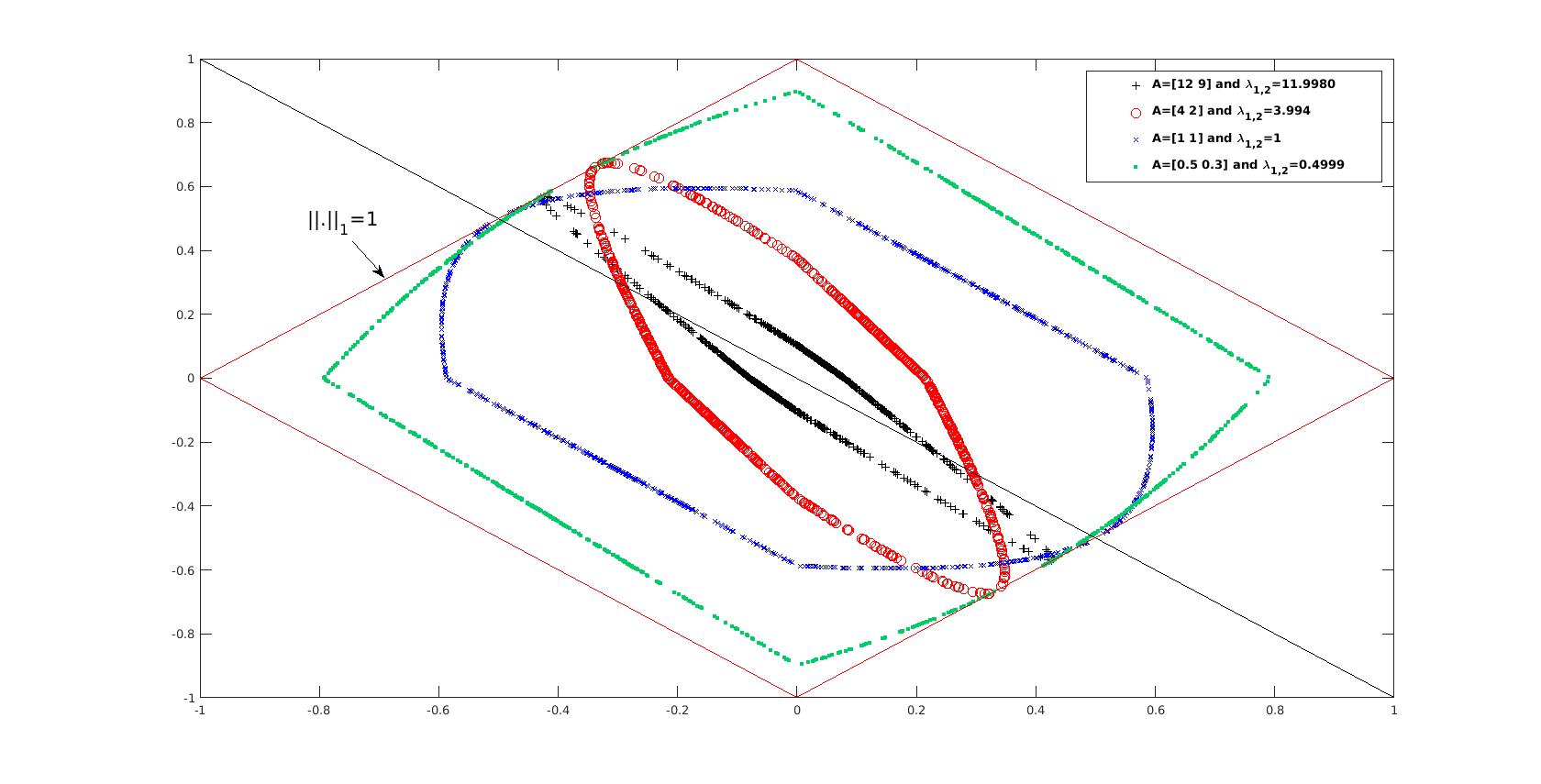} 
\caption{Contours of $\mathcal{B}(a_1,a_2,0)$ for different matrices $\A$, $p=2$ and $n=1$.} \label{Fig2}
\end{figure} 

\begin{remark} The numerics show that$Z(\omega_{lasso})$ exploses for the large values of $\omega_{lasso}$, it means that $\omega$ is closes to the null-space of $\A$. to eleminate that explosion we need to estimate the tail of the gaussian density.
Using the Gordon estimation \cite{Gordon^4}
\be 
\frac{\exp(-\frac{x^2}{2})}{x+\frac{1}{x}}\leq \sqrt{2\pi}(1-F(x))\leq 
\frac{\exp(-\frac{x^2}{2})}{x}, \quad x >0,
\ee 
we have the following approximation   
\ben\label{psivoisinagezero} 
\frac{1}{b^2}-\frac{1}{b}\leq z_2(b^2) \leq \frac{1}{b^2}-\frac{1}{1+b^3}, \quad near\,\, to\,\, 0.
\een 
\end{remark}
\section{MCMC diagnosis}
Here we take $p=7$, $n=4$, $\A \sim \Bc(\pm \frac{1}{\sqrt{n}})$ and 
for simplicity we consider $\y=0$. 
We sample from the distribution $c$ (\ref{c})
using Hastings-Metropolis algorithm  $(\x^{(t)})$ and propose the test 
$\|\x^{(t)}\|_2\leq q r(\theta^{(t)})$ 
as a criterion for the convergence. Here 
$\theta^{(t)}:=\frac{\x^{(t)}}{\|\x^{(t)}\|_2}$. 
We recall that if $\x$ is drawn from the target distribution $c$, then 
$\|\x\|_2\leq q r(\theta)$ with the probability at least equal to 
$P(q,p)$. Table 2 gives the values of the probability $P(q, p)$. Note that for $q \geq 2.5$ the criterion $\|\x^{(t)}\|_2\leq q r(\theta^{(t)})$ is satisfied with a large probability.

\renewcommand{\arraystretch}{0.9} 
\setlength{\tabcolsep}{0.07cm} 
 \begin{table}[!ht]
\begin{center}
\begin{tabular}{|c|c|c|c|c|c|c|c|}
\hline $q$ &2 &
 2.5 &
 3 &
   3.5 &
   4 &
   4.5&
   5\\
\hline $P(q,p)$ & 0.6672&
0.9446  &
0.9924 &
 0.9991&
 0.9999  &
 1.0000 &
 1.0000  \\
\hline  
\end{tabular}
\caption{ Values of the probability $P(q, p)$ for $p=7$.}
\end{center}
\end{table}

\subsection{Independent sampler (IS)}  
The proposal distribution 
\be 
Q(\x_2,\x_1)=p(\x_2)=\frac{1}{2^p}\exp(-\|\x_2\|_1),\quad \forall\,\x_1, \x_2.
\ee 
The ratio   
\be 
\frac{c(\x)}{p(\x)}\leq \frac{2^p}{Z},\quad \forall\,\x.
\ee 
It's known that MCMC $(x^{(t)})$ with the target distribution $c$ and the proposal distribution $p$ 
is uniformly ergodic \cite{Mengersen}:  
\be 
\sup_{A\subset \mathcal{B}(\Rb^p)}|\Pb(\x^{(t)}\in A\,|\,\x^{(0)})-\int_{A}c(\x)d\x|\leq
(1-\frac{Z}{2^p})^t. 
\ee
Here $Z\approx 2.2142$ and then $(1-\frac{Z}{2^p})=0.9827$. 
Figure 4(a) shows respectively the plot of $t\to 5r(\theta^{(t)})$ and $t\to \|\x^{(t)}\|_2$. 

\subsection{Random-walk (RW) Metropolis algorithm} 
We do not know if the target distribution $c$
satisfies the curvature condition in \cite{Roberts} Section 6. 
Here we propose  
to analyse the convergence of the Random walk Metropolis algorithm 
$(\x^{(t)})$ using the criterion $\|\x^{(t)}\|_2\leq q r(\theta^{(t)})$. 
Figure 4(b) shows respectively the plot of $t\to 5 r(\theta^{(t)})$ and $t\to \|\x^{(t)}\|_2$.  

Figures 4 show that contrary to independent sampler algorithm, 
the random walk (RW) algorithm satisfies early  
the criterion $\|\x^{(t)}\|_2\leq 5 r(\theta)$. More precisely   
\begin{itemize}
\item[1)] the independent sampler (IS) algorithm 
begins to satisfy the criterion $\|\x^{(t)}\|_2\leq 5 r(\theta^{(t)})$
at $t=8\times10^5$ iteration.

\item[2)] The RW algorithm begins to satisfy the criterion 
$\|\x^{(t)}\|_2\leq 3.5 r(\theta^{(t)})$
at $t=939065$ iteration, but the IS algorithm 
never satisfies the criterion $\|\x^{(t)}\|_2\leq 3.5 r(\theta^{(t)})$. 
\end{itemize}
We finally compare IS and RW algorithms using 
the fact that $\int_{\R^p} \x c(\x)d\x=0$. 
The best algorithm will furnish 
the best approximation of the integral $\int_{\R^p} \x c(\x)d\x$.
Table 3 gives the estimators 
$\frac{1}{N} \sum_{t=1}^{N} x_{IS}^{(t)}\approx \int_{\R^p} \x c(\x)d\x$ and 
$\frac{1}{N} \sum_{t=1}^{N} x_{RW}^{(t)}\approx \int_{\R^p} \x c(\x)d\x$. 
It follows that $\|\frac{1}{N} \sum_{t=1}^{N} x_{IS}^{(t)}\|_2 = 0.0187$ and $\|\frac{1}{N} \sum_{t=1}^{N} x_{RW}^{(t)}\|_2 = 0.0041$. We conclude that the random walk algorithm wins 
for both criteria against independent sampler algorithm.   

\renewcommand{\arraystretch}{0.9} 
\setlength{\tabcolsep}{0.07cm} 
 \begin{table}[!ht]
\begin{center}
\begin{tabular}{|c|c|c|c|c|c|c|c|}
\hline  &$x_1$ &
 $ x_2$ &
   $x_3$&
   $x_4$&
   $x_5$&
   $x_6$&
   $x_7$\\
\hline $\hat{x}_{IS}$ &-0.0005&
  -0.0037&
 0.0016&
 0.0164&
   0.0050&
  0.0021&
  -0.0058 \\
\hline  
 $\hat{x}_{RW}$&0.0005 &
-0.0019&
 -0.0002&
 0.0012&
  -0.0005&
 0.0031&
  -0.0011 \\
\hline 
\end{tabular}
\caption{ $IS$ and $RW$ estimators using $N=10^6$ iterations.}
\end{center}
\end{table}

\begin{figure}[!ht]\centering
\includegraphics[width=16cm]{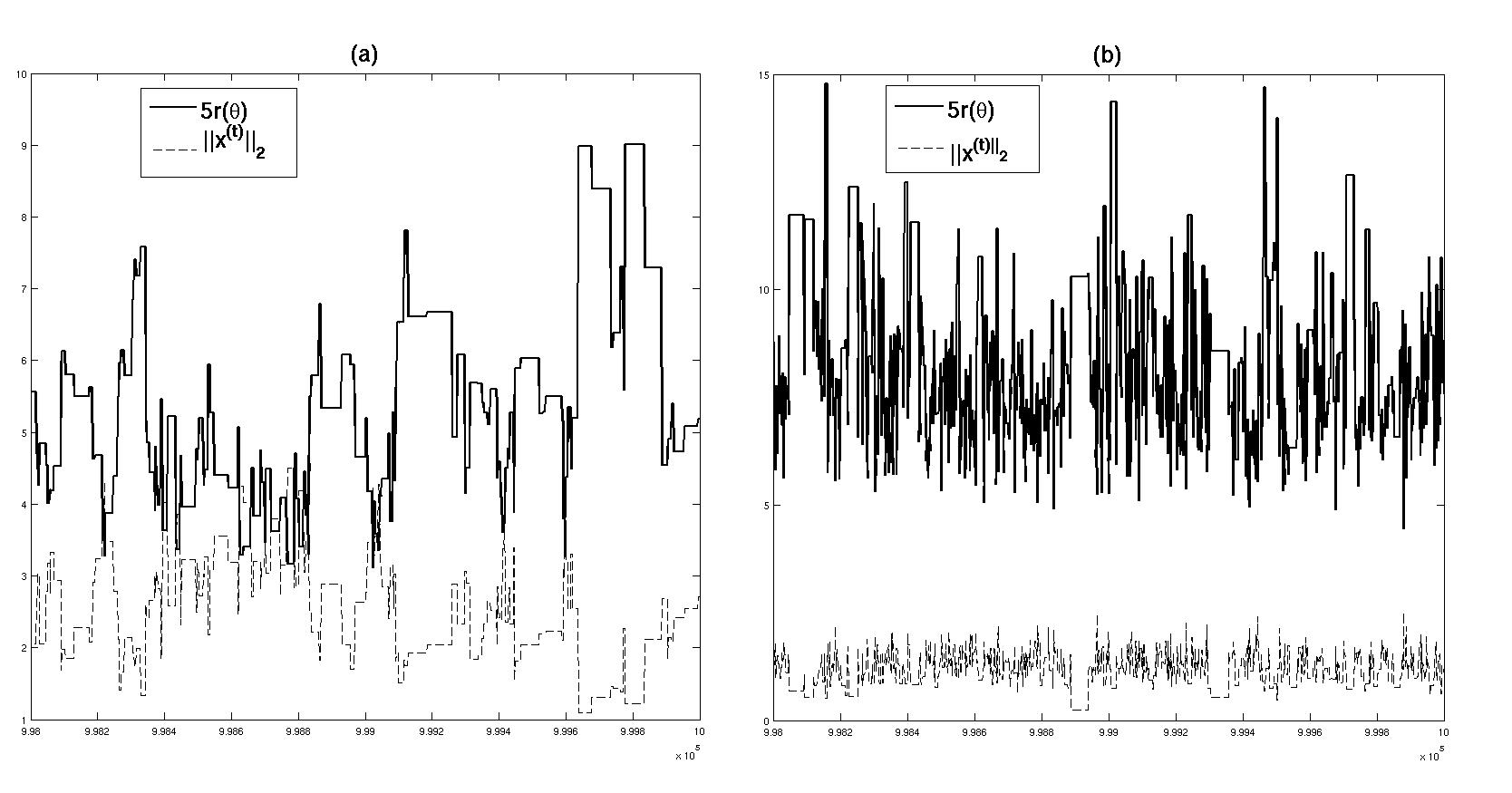} 
\caption{(a): Test of convergence of MCMC algorithm with proposal distribution $p(\x_2)$. (b): Test of convergence of MCMC algorithm with $\Nc(0, 0.5\I_p)$ proposal distribution. $N=10^6$ iterations.}
\end{figure}    

\section{Conclusion}
We studied the geometry of bayesian LASSO using polar coordinates  
and calculated the partition function. 
We obtained a concentration inequality and derived MCMC convergence diagnosis
for the convergence of Hasting Metropolis algorithm. 
We showed that the random walk MCMC with the variance 0.5 wins again 
the independent sampler with Laplace proposal distribution.

\end{document}